\documentclass[12pt]{article}
\usepackage{graphicx}
\usepackage{amssymb}
\usepackage{amsmath}
\usepackage{amsthm}
\usepackage{mathptmx}
 \usepackage{latexsym}

\usepackage{tikz}
\usetikzlibrary{decorations.pathreplacing,calc}

\usepackage{here}
\usepackage[]{hyperref}

\def\<#1>{\langle#1\rangle}

\newtheorem{thm}{Theorem}[section]
  \newtheorem{definition}[thm]{Definition}

   \theoremstyle{definition}

\begin{document}

\title{The arithmetic of simplices}

\author{Edward Mieczkowski\footnote{
 Gdynia, Poland;  E-mail: \texttt{edmieczyk@gmail.com}.}}

\date{}

\maketitle

\begin{abstract}
\noindent This paper continues the study initiated in "The aithmetic of Triangles." \
We begin by examining a set of similar tetrahedra with parallel sides, together with a set of points in three-dimensional space. It turns out that the set

\noindent $\mathbb{R}_3= \{\pm  \<x >=\pm (x^3,x^2,x,1); x\in\mathbb{R} \}$\\
effectively characterizes this family of tetrahedra. The set $\mathbb{R}_3$ is a subset of the ring\\
 $\mathbb{R}^4 = \mathbb{R} \times \mathbb{R} \times \mathbb{R} \times \mathbb{R} = { (x, y, z, w) ; ; ; x, y, z, w \in \mathbb{R} }$,\\
  with addition and multiplication defined component-wise.
The set $\mathbb{R}_3$ supports two operations. Multiplication is inherited directly from the ring $\mathbb{R}^4$, while addition is a four-argument  operation that reflects geometric transformations such as homothety and translation of elements in $\mathbb{R}_3$. Depending on the relative magnitudes of the arguments in this addition, there are 15 distinct geometric interpretations.

\noindent However, the defined addition has its limitations. It turns out that, within this framework, the reduction of terms with different signs is not always possible. This leads to the distinction between an equation that is true in the arithmetic sense and one that is true in the geometric sense.

\noindent A novel form of addition in $\mathbb{R}_3$ leads to intriguing properties of multiplication in $\mathbb{R}_3$, which are examined in a dedicated chapter.
 
\noindent We then generalize this approach to sets of $k$-dimensional similar simplices with parallel sides, along with corresponding sets of points in $k$-dimensional space. These can be described by the set\\
$\mathbb{R}_k= \{\pm  \<x >=\pm (x^k,x^{k-1},\ldots,x,1); x\in\mathbb{R} \}$ \\
a subset of the ring \\
$\mathbb{R}^{k+1}
= \{ (x_1,x_2,\ldots,x_{k+1}) ; \; x_1,x_2,\ldots,x_{k+1}\in\mathbb{R} \}$\\
with addition and multiplication again defined component-wise.
\end{abstract}

\newpage
\subsection*{ Introduction}

Let us take the ring $\mathbb{R}^4=\mathbb{R}\times\mathbb{R}\times\mathbb{R}\times\mathbb{R}
= \{ (x,y,z,w) ; \; x,y,z,w\in\mathbb{R} \}$ with  additon and multiplication\\
\begin{align}
(x_1,y_1,z_1,w_1)+(x_2,y_2,z_2,w_2)=&\ (x_1+x_2,y_1+y_2,z_1+z_2,w_1+w_2), \label{1} \\
(x_1,y_1,z_1,w_1)\cdot (x_2,y_2,z_2,w_2)=&\ (x_1\cdot x_2,y_1\cdot y_2,z_1\cdot z_2,w_1\cdot w_2). \label{2}
\end{align}
Let us consider the subset of the ring $\mathbb{R}^4$,\\
 the set $\mathbb{R}_3= \{\pm  \<x >=\pm (x^3,x^2,x,1); x\in\mathbb{R} \}$.\\
It is closed under multiplication (\ref{2}) but not under addition (\ref{1}).\\
The set $\mathbb{R}_3$ is closed under    the following kind of addition
(It should be noted that the summands in below Eq.(\ref{3}) are not  the coordinates of the ring $\mathbb{R}^{4}$).
   
 \begin{align}\begin{split}
\forall x,y,z,w,t \in \mathbb{R}\qquad & \\ 
       \<x+y+z+w+t> = &\ \< x\!+\!y\!+\!z\!+\!t > + \< x\!+\!y\!+\!w\!+\!t> + \< x\!+\!z\!+\!w\!+\!t >
       					 +\< y\!+\!z\!+\!w\!+\!t > \\
       				&\ - \< x+y+t > - \< x+z+t> - \< x+w+t > \\
                    &-\< y+z+t >-\< y+w+t > -\< z+w+t > \\
                    & + \< x+t> + \<y+t> + \<z+t> +\< w+t >-\<t>   \label{3}
\end{split}\end{align}
because the equations
\begin{align*}
\forall x,y,z,t \in \mathbb{R} \qquad \; \; & \ \forall i=3,2,1,0\qquad \qquad \\
       (x\!+\!y\!+\!z\!+\!w\!+\!t)^i = &\ ( x\!+\!y\!+\!z\!+\!t )^i + 
       ( x\!+\!y\!+\!w\!+\!t)^i + ( x\!+\!z\!+\!w\!+\!t)^i + (y\!+\!z\!+\!w\!+\!t )^i \\
       		      & - ( x+y+t )^i - ( x+z+t)^i - ( x+w+t )^i\\
       		      & - ( y+z+t )^i -( y+w+t )^i-( z+w+t )^i  \\
                    &\ + ( x+t)^i + (y+t)^i + (z+t)^i + ( w+t )^i-(t)^i
   \end{align*}
are true.\\
If we multiply the Eq. (\ref{3}) by $-1$, we get the definition of addition for the elements $-\<x>$ of the set $\mathbb{R}_3$.

\noindent It can be proved by mathematical induction that the following equation follows   from   Eq. (\ref{3}).  
 \begin{align}
\forall n \in \mathbb{N} \quad 
   \<n>= &\ \frac{(n-1)n(n+1)}{6}\<2>-\frac{(n-2)n(n+1)}{2}\<1> \nonumber\\
        &\ + \frac{(n-2)(n-1)(n+1)}{2}\<0>	-\frac{(n-2)(n-1)n}{6}\<-1>.
    \label{4}
\end{align}

\noindent But you can check purely for accounting that  Eq. (\ref{4}) holds for each $x \in \mathbb{R}$.
\begin{align}
\forall x \in \mathbb{R}\quad
       \<x>=&\ \frac{(x-1)x(x+1)}{6}\<2>-\frac{(x-2)x(x+1)}{2}\<1> \nonumber\\
        &\ + \frac{(x-2)(x-1)(x+1)}{2}\<0>	-\frac{(x-2)(x-1)x}{6}\<-1>.
    \label{5}
\end{align}
Let us transform Eq. (\ref{5}). 
\begin{align*}
\< x > = & \ \frac{x^3-x}{6}\< 2> -\frac{x^3-x^2-2x}{2}\<1> + \frac{x^3-2x^2-x+2}{2}\<0>
		-\frac{x^3-3x^2+2x}{6}\<-1> \\
		= &\ x^3\frac{\<2>-3\< 1> +3\<0>-\< -1>}{6}+  x^2\frac{\< 1> -2\<0>+\< -1>}{2} \\
		  &\ +x\frac{-\<2>+6\< 1> -3\<0> -2\< -1>}{6} + 1\<0>.
\end{align*}
 It is easy to check that 
 the elements \\[.1cm] 
   $\frac{\<2>-3\< 1> +3\<0>-\< -1>}{6} = A_3$, $\frac{\< 1> -2\<0>+\< -1>}{2} = A_2$, 
  $\frac{-\<2>+6\< 1> -3\<0> -2\< -1>}{6} = A_1$, $\<0>=A_0$ \\[.1cm]
  are orthogonal. So we can put\\[.1cm]
$ A_3=(1,0,0,0)$,\; $ A_2=(0,1,0,0)$, \; $A_1=(0,0,1,0)$, \; $A_0=(0,0,0,1)$ and now we know why $\<x>=(x^3,x^2,x,1)$.\\

\section{Geometric interpretation of the set $\mathbb{R}_3$.} 

\noindent  Let us set any closed  tetrahedron in the 3-dimensional  space (the tetrahedron does not have to be regular) and let us denote it by the symbol $\< 1>$. Each tetrahedron obtained from the fixed tetrahedron $\<1>$ by any translation will be denoted by $\<1>$ too.\\
The tetrahedron obtained from the tetrahedron $\<1>$ by any homothety of ratio $x>0$ \cite{Wol} will be denoted by $\<x>$.\\
Each point of the space will be denoted by $\<0>$.\\
The symbol  $-\< x>$ denotes the tetrahedron, which if is a cover of  the tetrahedron $\< x>$ in topological sense  \cite{Wol} and vice versa  the tetrahedron $\< x>$ is a cover of  the tetrahedron $-\< x>$ then both give an empty set  denoted by an element $(0,0,0,0) \in \mathbb{R}^4$.\\
We will write easier that if the tetrahedron $-\< x>$ lies on  the tetrahedron $\< x>$ then both give an empty set.
We will  call the tetrahedron  $-\< x>$ sometimes antytetrahedron.\\
An element $\<x>$ for each $x \in \mathbb{R},\: x>0$  will be interpreted as a closed tetrahedron and marked with a black color. 
\tikz{\draw[fill,gray!50] (0,0)--(1,0)--(.5,.8)--(0,0);
\draw[dashed](0,0)--(.5,.25)--(1,0);\draw[dashed](.5,.25)--(.5,.8);
\draw (0,0)--(1,0)--(.5,.8)--cycle;
\draw[fill](0,0)circle(1pt);\draw[fill](1,0)circle(1pt);
\draw[fill](.5,.8)circle(1pt);\draw[fill](.5,.25)circle(1pt);} \ .
 The interior of this
tetrahedron should be black but we would like to mark that edges and vertexes belong to  this tetrahedron,  so we coloured  it in gray.\\
An element $-\<x>$ will be interpreted as a closed tetrahedron  and marked with a red color
\tikz{\draw[fill,red!20] (0,0)--(1,0)--(.5,.8)--(0,0);
\draw[dashed,red](0,0)--(.5,.25)--(1,0);\draw[dashed,red](.5,.25)--(.5,.8);
\draw[red] (0,0)--(1,0)--(.5,.8)--cycle;
\draw[fill,red](0,0)circle(1pt);\draw[fill,red](1,0)circle(1pt);
\draw[fill,red](.5,.8)circle(1pt); \draw[fill,red](.5,.25)circle(1pt);} \;.\\
An element $\<0>$ will be marked as a black point and an element $-\<0>$ will be marked as a red point.\\
To see the result of putting the black tetrahedron on the red one, we will denote the empty set of green.\\
An element $\<-x>$ for each $x \in \mathbb{R},\: x>0$  will be interpreted as a open tetrahedron (that is, without faces, edges and vertices) symmetrical to the  tetrahedron  $\<x >$ in relation to any edge and will be marked with a red color inside and with green faces, edges, and vertices.
\tikz{\draw[fill,red!20] (.5,0)--(1,.8)--(0,.8)--(.5,0);
\draw[green] (.5,0)--(1,.8)--(.5,.55)--(.5,0)--(0,.8)--(.5,.55)--(1,.8)--(0,.8);
\draw[fill,green](.5,0)circle(1pt);\draw[fill,green](1,.8)circle(1pt);
\draw[fill,green](0,.8)circle(1pt);\draw[fill,green](.5,.55)circle(1pt);} .\\
At the end an element $-\<-x>$ for each $x \in \mathbb{R},\: x>0$  will be interpreted as an open tetrahedron symmetrical to the  tetrahedron  $-\<x >$ in relation to any edge and will be marked with a black color inside and with green faces, edges, and vertices.
\tikz{\draw[fill,gray!50] (.5,0)--(1,.8)--(0,.8)--(.5,0);
\draw[green] (.5,0)--(1,.8)--(.5,.55)--(.5,0)--(0,.8)--(.5,.55)--(1,.8)--(0,.8);
\draw[fill,green](.5,0)circle(1pt);\draw[fill,green](1,.8)circle(1pt);
\draw[fill,green](0,.8)circle(1pt);\draw[fill,green](.5,.55)circle(1pt);} .\\

\noindent Why the tetrahedron $\<-x>$ and $-\<-x>$ for $ x>0$  are open we will see after introducing the geometric interpretation of the operation (\ref{3}).\\

\noindent The geometric interpretation of the individual components of Eq. (\ref{4}) is as follows (Fig. \ref{f1}).
 
\begin{figure}[H]
\begin{center} 
\begin{tikzpicture}[>=stealth,scale=.5]
\draw[fill,gray!30] (0,0)--(8,16)--(16,0)--(0,0);
\draw(0,0)--(8,16)--(16,0)--(0,0);
\draw[dashed] (0,0)--(8,4)--(16,0);
\draw[dashed] (8,4)--(8,16);
\draw[dotted] (2,1)--(14,1)--(12,0)--(6,3)--(10,3)--(4,0)--cycle;
\draw[dotted] (4,2)--(8,0)--(12,2)--cycle;
\draw[dotted](2,4)--(2,1)--(8,13)--(6,12)--(6,3)--(8,7)--cycle;
\draw[dotted] (4,8)--(4,2)--(8,10)--cycle;
\draw[dotted] (14,1)--(14,4)--(8,7)--(10,3)--(10,12)--(8,13)--cycle;
\draw[dotted] (12,2)--(12,8)--(8,10)--cycle;
\draw[dotted,thick] (2,4)--(4,0)--(10,12)--(6,12)--(12,0)--(14,4)--cycle;
\draw[dotted,thick] (4,8)--(8,0)--(12,8)--cycle;
\draw[dotted,red,thick] (4,5)--(12,5)--(10,4)--(6,6)--(10,6)--(6,4)--cycle;
\draw[dotted,red,thick] (12,5)--(10,1)--(6,9)--(10,9)--(6,1)--(4,5);
\draw[dotted,red,thick] (6,9)--(8,8)--(8,2)--(6,6)--(10,4)--(10,1);
\draw[dotted,red,thick] (10,6)--(8,2)--(8,8)--(10,9)--(6,1)--(6,4);

\draw [decorate,decoration={brace,amplitude=9pt},
       xshift=6pt,yshift=2pt]
   (8,16)--(12,8); 
\draw node[] at (11.5,12.5){\small{$\<2>$}};   
\draw node[] at (7.5,13.5){\small{$\<1>$}};
\draw node[] at (10.5,9.5){\small{$\<1>$}};
\draw node[] at (1.5,1.5){\small{$\<1>$}};
\draw node[] at (10.8,3.7){\scriptsize{$\text{-}\<\text{-}1>$}};
\draw node[] at (8.8,7.7){\scriptsize{$\text{-}\<\text{-}1>$}};
 \end{tikzpicture}
\caption{} \label{f1}
\end{center}
\end{figure}
\noindent If the  tetrahedron  $\<1 >$  has a volume of 1, then the  tetrahedron  $\<n >$  has a volume of $n^3$.\\
 Each face of the tetrahedron $\<n>$ has $n^2$  of all tetrahedra $\<1>$ and $\<-1>$.\\
 If we assume that the length of each edge of the tetrahedron $\< 1> $ is equal to $1$ (even when the triangle is not regular) then $n$ is the length of each edge of the tetrahedron  $\<n>$.\\
The number 1 in $(n^3,n^2,n,1)$ means one  tetrahedron.\\

\noindent\textit{Geometric interpretation of the operation (\ref{3}).}\\ 

\noindent Geometric interpretation of the adding (\ref{3}) is similar to adding in the set $\mathbb{R}_2$  
 \cite{Art}.\\ 
 \noindent Let us fix the ordered successive components   $ x,y,z,w $ of the sum \mbox{$ \< x+y+z+w+ t > $.} This components extend or shorten the tetrahedron $ \< t> $ in the directions \mbox{I, II, III, IV} or \mbox{I ', II ', III', IV'}, depending on whether the numbers $ x, y, z, w $     are positive or negative (Fig.~\ref{f2}).
 
\begin{figure}[H]
\begin{center} 
\begin{tikzpicture}[>=stealth,scale=1.2]
\draw[fill,gray!30] (0,0)--(3,6)--(6,.5)--(0,0);
\draw(0,0)--(3,6)--(6,.5)--(0,0);
\draw[dashed] (0,0)--(2,1.5)--(6,.5);
\draw[dashed] (2,1.5)--(3,6);

\draw[thick,->](3.5,2.7)--(4,3);
\draw[fill](3.5,2.7)circle(1pt);
\draw[thick,dotted,->](3.5,2.7)--(3,2.4);
\draw node[] at (4.1,3){\small{I}};
\draw node[] at (2.9,2.4){\small{I'}};
\draw[thick,->](1.8,2.7)--(1.3,3);
\draw[fill](1.8,2.7)circle(1pt);
\draw[thick,dotted,->](1.8,2.7)--(2.3,2.4);
\draw node[] at (1.1,3.1){\small{II}};
\draw node[] at (2.5,2.3){\small{II'}};
\draw[thick,->](2.8,1.7)--(3.1,1.5);
\draw[fill](2.8,1.7)circle(1pt);
\draw[thick,dotted,->](2.8,1.7)--(2.5,1.9);
\draw node[] at (3.3,1.4){\small{III}};
\draw node[] at (2.3,2){\small{III'}};
\draw[thick,->](2.2,.7)--(2.2,-.1);
\draw[fill](2.2,.7)circle(1pt);
\draw[thick,dotted,->](2.2,.7)--(2.2,1.1);
\draw node[] at (2.2,-.3){\small{IV}};
\draw node[] at (2.2,1.25){\small{IV'}};
 \end{tikzpicture}
\caption{} \label{f2}
\end{center}
\end{figure} 

\noindent Therefore Eq. (\ref {3}) should be properly written as
\begin{align*}\begin{split}
\forall x,y,z,w,t \in \mathbb{R}\qquad & \\ 
    \<x+y+z+w+t> = &\ \< x\!+\!y\!+\!z\!+\!0\!+\!t > + \< x\!+\!y\!+\!0\!+\!w\!+\!t> \\
   					&\	 + \< x\!+\!0\!+\!z\!+\!w\!+\!t > +\<0\!+\!y\!+\!z\!+\!w\!+\!t > \\
       	&\ -\< x\!+\!y\!+\!0\!+\!0\!+\!t > -\< x\!+\!0\!+\!z\!+\!0\!+\!t >-\< x\!+\!0\!+\!0\!+\!w\!+\!t >
       				 \\
        &\ - \< 0\!+\!y\!+\!z\!+\!0\!+\!t >-\< 0\!+\!y\!+\!0\!+\!w\!+\!t >-\< 0\!+\!0\!+\!z\!+\!w\!+\!t >
        			 \\
        &\ + \< x\!+\!0\!+\!0\!+\!0\!+\!t >+\< 0\!+\!y\!+\!0\!+\!0\!+\!t >  \\
        &\ + \< 0\!+\!0\!+\!z\!+\!0\!+\!t >+\< 0\!+\!0\!+\!0\!+\!w\!+\!t >-\< 0\!+\!0\!+\!0\!+\!0\!+\!t >          
        \end{split}   \end{align*}

\noindent Since the above record is long and inconvenient, we will replace it with the following one

\begin{align}\begin{split}
  \forall x,y,z,w,t \in \mathbb{R} \qquad & \\
   \<x,y,z,w,t>&= \< x,y,z,0,t> + \<x,y,0,w,t> + \< x,0,z,w,t>+\<0,y,z,w,t>  \\
             &\quad - \<x,y,0,0,t> -\< x,0,z,0,t> - \< x,0,0,w,t > \\
             &\quad - \<0,y,z,0,t>- \<0,y,0,w,t>- \<0,0,z,w,t>\\
             &\quad + \<x,0,0,0,t>+ \<0,y,0,0,t>+ \<0,0,z,0,t>+ \\
             &\quad +\<0,0,0,w,t> - \< 0,0,0,0,t >,
    \label{6}
 \end{split}   \end{align}
 where we can write $\< 0,0,0,t >$   as $\< t >$.\\
 From now on, we will replace Eq. (\ref{3}) by Eq. (\ref{6}).\\
\noindent Below we have  examples of the creation of new tetrahedra from the tetrahedron $ \< t >, \; t>0 $ (Fig.~\ref{f3}-\ref{f6}).\\
Because some components of Eq. (\ref{6}) are positive (i.e. black) and others are negative (i.e. red),  for simplicity we  omit these colors in Fig. ~\ref{f3}-\ref{f6}.\\
The tetrahedron $ \<t> $ has vertices $ ABCD $, and  the tetrahedron $ \<x,y,z,t>$  has vertices $ A'B'C'D'$.\\ 
The faces created by moving the faces $BCD$, $ACD$, $ABD$ and $ABC$ by values $x,\ y, \ z$,  and $w$ are blue, yellow, orange and brown respectively.

\begin{figure}[H]
\begin{minipage}[t]{.45\textwidth}
 \begin{center}
\begin{tikzpicture}[>=stealth,scale=.5]

\draw(3,6)--(4.2,8.4);\draw[dashed](2,1.5)--(2.8,2.1);
\draw(6,.5)--(8.4,.7);
\draw(0,0)--(3,6)--(6,.5)--(0,0);
\draw[blue,thick,dashed](4.2,8.4)--(2.8,2.1)--(8.4,.7);
\draw[blue,thick](4.2,8.4)--(8.4,.7);
\draw[dashed] (0,0)--(2,1.5)--(6,.5);
\draw[dashed] (2,1.5)--(3,6);
\draw node at (-.2,-.4){\footnotesize $A'=A$};\draw node at (6,.1){\footnotesize $B$};
\draw node at (8.4,.3){\footnotesize $B'$};
\draw node at (1.7,1.75){\footnotesize $C$};\draw node at (3.3,2.4){\footnotesize $C'$};
\draw node at (2.7,6.3){\footnotesize $D$};\draw node at (4.1,8.7){\footnotesize $D'$};
\end{tikzpicture}
\caption{The tetrahedron $ A'B'C'D' = \< x,0,0,0,t>$ for  $x>0$.} \label{f3}
\end{center}
\end{minipage}
\hfill
\begin{minipage}[t]{.35\textwidth}
 \begin{center}
 \begin{tikzpicture}[>=stealth,scale=0.5]
\draw(0,0)--(3,6)--(6,.5)--(0,0);
\draw[yellow,thick,dashed](2.4,.2)--(3.8,1.1)--(4.2,3.8);
\draw[yellow,thick](2.4,.2)--(4.2,3.8);
\draw[dashed] (0,0)--(2,1.5)--(6,.5);
\draw[dashed] (2,1.5)--(3,6); 
\draw node at (-.1,-.4){\footnotesize $A$};\draw node at (2.4,-.2){\footnotesize $A'$};
\draw node at (5.8,.1){\footnotesize $B=B'$};
\draw node at (1.7,1.75){\footnotesize $C$};\draw node at (4.3,1.35){\footnotesize $C'$};
\draw node at (2.7,6.3){\footnotesize $D$};\draw node at (4.6,4.2){\footnotesize $D'$};
\end{tikzpicture}
\caption{The tetrahedron $ A'B'C'D' = \< 0,y,0,0,t>$ for  $-t<y<0$.} \label{f4}
 \end{center}
\end{minipage}
\end{figure} 

\begin{figure}[H]
\begin{minipage}[t]{.33\textwidth}
 \begin{center}
\begin{tikzpicture}[>=stealth,scale=.5]
\draw(0,0)--(3,6)--(6,.5)--(0,0);
\draw[dashed] (0,0)--(2,1.5)--(6,.5);
\draw[dashed] (2,1.5)--(3,6);
\draw(3,6)--(4.2,8.4);\draw(3,6)--(1.8,8.2);
\draw(3,6)--(3.4,7.8);
\draw[brown,thick](4.2,8.4)--(1.8,8.2)--(3.4,7.8)--cycle;
\draw node at (0,-.4){\footnotesize $A$};\draw node at (6,.1){\footnotesize $B$};
\draw node at (1.5,8.3){\footnotesize $B'$};
\draw node at (1.7,1.75){\footnotesize $C$};\draw node at (3,7.5){\footnotesize $C'$};
\draw node at (4.2,6.05){\footnotesize $D=D'$};\draw node at (4.5,8.8){\footnotesize $A'$};
\end{tikzpicture}
\caption{The tetrahedron $ A'B'C'D' = \< 0,0,0,w,t>$ for  $w<-t$.} \label{f5}
\end{center}
\end{minipage}
\hfill
\begin{minipage}[t]{.48\textwidth}
 \begin{center}
 \begin{tikzpicture}[>=stealth,scale=0.5]
\draw(0,0)--(3,6)--(6,.5)--(0,0);
\draw[dashed] (0,0)--(2,1.5)--(6,.5);
\draw[dashed] (2,1.5)--(3,6); 
\draw[dashed](0,0)--(-.8,-.6);\draw[dashed](3,6)--(3.4,7.8);
\draw[dashed](6,.5)--(7.6,.1);
\draw[orange,thick](-.8,-.6)--(3.4,7.8)--(7.6,.1)--cycle;
\draw[orange,thick](-.8,-.6)--(-3.2,-.8);
\draw[orange,thick,dashed](-3.2,-.8)--(-.8,-5.2);
\draw[blue,thick,dashed](-.8,-5.2)--(-3.2,-.8);
\draw[orange,thick,dashed](-2,-3)--(-.8,-.6);
\draw[yellow,thick,dashed](2.4,.2)--(3.6,1.1)--(4.2,3.8);
\draw[yellow,thick,dashed](4.2,3.8)--(4.6,5.6)--(-.8,-5.2);
\draw[orange,thick,dashed](-.8,-5.2)--(4.6,5.6);

\draw[yellow,thick](2.4,.2)--(4.2,3.8);
\draw(0,0)--(-4.8,-.4);\draw[dashed](0,0)--(-2.4,-4.8);
\draw[](-.8,-.6)--(-1.6,-1.2);
\draw[blue,thick](-4.8,-.4)--(-2.4,-4.8)--(-1.6,-1.2)--cycle;
\draw[blue,thick](-1.6,-1.2)--(0,-1.6);
\draw[blue,thick](-2.4,-4.8)--(-1.2,-7);
\draw[blue,thick,dashed](-1.2,-7)--(0,-1.6);
\draw[yellow,thick,dashed](0,-1.6)--(-1.2,-7);
\draw[yellow,thick](2.4,.2)--(0,-1.6);\draw[yellow,thick](2.4,.2)--(-1.2,-7);
\draw node at (-.1,.4){\footnotesize $A$};\draw node at (1.7,-.7){\footnotesize $A'$};
\draw node at (6.3,.75){\footnotesize $B$};\draw node at (-3.5,-1.2){\footnotesize $B'$};
\draw node at (1.7,1.75){\footnotesize $C$};\draw node at (.35,-1.8){\footnotesize $C'$};
\draw node at (3.4,6.1){\footnotesize $D$};\draw node at (-1.3,-5.4){\footnotesize $D'$};
\end{tikzpicture}
\caption{The tetrahedron $ A'B'C'D' = \< x,y,z,0,t>$ for $x<-t$, $-t<y<0$, $z>0$.} \label{f6}
 \end{center}
\end{minipage}
\end{figure} 

\noindent From the below equation
\[
\< -1>=\<-1,-1,-1,-1,3>=4\<0>-6\<1>+4\<2>-\<3>
\]
we  can get the way  of building of the triangle $\<-1>$.
Fig. \ref{f7} shows each step of  the construction of $\<-1>$.\\

\noindent Similarly,  from the relation
\[ \<-x>=\<-x,-x,-x,-x,3x>= 4\<0>-6\<x>+4\<2x>-\<3x>\]
we can see that for each $x \in \mathbb{R},\; x>0$ the tetrahedron  $ \< -x>$   is the opened tetrahedron.\\

\begin{figure}[H]
\begin{center} 
\begin{tikzpicture}[>=stealth,scale=.8]
\draw[fill,red!20] (0,0)--(3,6)--(6,0)--cycle;
\draw[red](0,0)--(3,6)--(6,0)--cycle;
\draw[red,dashed] (0,0)--(3,1.5)--(6,0);
\draw[red,dashed] (3,1.5)--(3,6);
\draw node[] at (3,-.6){\small{$-\<3>$,}};

\draw[fill,red!20](8,2.5)--(10,2.5)--(9,.5)--cycle;
\draw[fill,gray!30](8,0)--(10,0)--(9,2)--cycle;
\draw[fill,gray!30](7,2)--(9,2)--(8,4)--cycle;
\draw[fill,gray!30](9,2)--(11,2)--(10,4)--cycle;
\draw[fill,gray!30](7,.5)--(9,.5)--(8,2.5)--cycle;
\draw[fill,gray!30](9,.5)--(11,.5)--(10,2.5)--cycle;
\draw[fill,gray!30](8,2.5)--(10,2.5)--(9,4.5)--cycle;
\draw[dashed](8,2.5)--(9,2)--(10,2.5);\draw[dashed](9,2)--(11,2)--(9,3)--(7,2)--cycle;
\draw[dotted](9,.5)--(10,2.5)--(8,2.5)--cycle;
\draw(8.75,2.5)--(9.25,2.5);
\draw(9.38,1.26)--(9.75,2);
\draw(8.62,1.26)--(8.25,2);
\draw(8.38,3.26)--(9,4.5)--(9.62,3.26);
\draw(8.26,.5)--(7,.5)--(7.75,2);
\draw(9.74,.5)--(11,.5)--(10.25,2);
\draw[dashed](9,.5)--(9,2)--(10,4)--(10,1)--(8,0)--(9,2);
\draw[dashed](9,.5)--(8,1)--(8,4)--(10,0)--(9,.5);
\draw[dotted](7,.5)--(11,.5)--(9,4.5)--cycle;
\draw[](7,2)--(11,2)--(10,4)--(8,0)--(10,0)--(8,4)--cycle;
\draw[dashed](9,4.5)--(9,3);\draw[dashed](7,2)--(8,2.5);
\draw[dashed](10,4)--(11,2);\draw[dashed](7,.5)--(8,1);
\draw[dashed](10,1)--(11,.5);
\draw[fill](9,2)circle(1pt);\draw[fill](8,2.5)circle(1pt);
\draw[fill](10,2.5)circle(1pt);\draw[fill](9,.5)circle(1pt);
\draw node[] at (9.3,1.8){\small{2}};\draw node[] at (10.2,2.65){\small{2}};
\draw node[] at (7.8,2.65){\small{2}};\draw node[] at (9,.25){\small{2}};
\draw node[] at (9,-.6){\small{$-\<3>+4\<2>$,}};

\draw[fill,red!20](12,2.5)--(14,2.5)--(13,.5)--cycle;
\draw[green,thick](12,2.5)--(14,2.5)--(13,.5)--cycle;
\draw[green,thick](12,2.5)--(13,2)--(13,.5);
\draw[green,thick](14,2.5)--(13,2)--cycle;
\draw[fill,red](13,2)circle(1pt);\draw[fill,red](12,2.5)circle(1pt);
\draw[fill,red](14,2.5)circle(1pt);\draw[fill,red](13,.5)circle(1pt);
\draw node[] at (13,-.6){\small{$-\<3>+4\<2>-6\<1>$,}};

\draw[fill,red!20](15,2.5)--(17,2.5)--(16,.5)--cycle;
\draw[green,thick](15,2.5)--(17,2.5)--(16,.5)--cycle;
\draw[green,thick](15,2.5)--(16,2)--(16,.5);
\draw[green,thick](17,2.5)--(16,2)--cycle;
\draw node[] at (16,-.6){\small{$\<-1>$}};
 \end{tikzpicture}
\caption{Stages of  the construction of the tetrahedron $\<-1>$.
 Numbers 2  mean an overlapping two points. The tetrahedron $\<-1>$ is red  inside and has green faces, edges, and vertices.} \label{f7}
\end{center}
\end{figure}

\section{The equation true in  the geometric sense.}

\noindent Let us consider Eq. (\ref{6}) for  concrete numbers
\begin{align}
  \<7> = \<3,2,1,1,0 > =&\ 2\<6>+\<5>+\<4> -\<5>-2\<4>-2\<3>-\<2> \nonumber \\
  					  &\ +\<3>+\< 2> + 2\<1> - \<0>. \label{7}
\end{align}
After reduction we will get
\begin{align}
    \<7> =  & 2\< 6> -\<4>-\<3>+2\<1>-\<0>. \label{8}
\end{align}
\noindent From the arithmetic point of view Eq.  (\ref{8})
 is true. But it easy to see that we can  not build the tetrahedron $\<7> $
 using only  tetrahedra from Eq. (\ref{8}). We can not even build the face of the tetrahedron $\<7> $. \\ 

\noindent So it makes sense to introduce the following definitions.
\begin{definition}
The equation $\< x> = \sum_j\alpha_j\<x_j>$, where $\alpha_j \in \{-1,1\}$ is true in the geometric sense if we can build the tetrahedron $\< x>$ from the elements $\alpha_j\< x_j>$. \label{d2.1}
\end{definition}

\begin{definition}
The equation $\< x> = \sum_j\alpha_j\<x_j>$, where $\alpha_j \in \{-1,1\}$ is true in the arithmetic sense if the equations  $ x^i  = \sum_j\alpha_j x^i_j$, where $\alpha_j \in \{-1,1\}$,  $i=0,1,2,3$ \label{d2.2}
hold.
\end{definition}

\noindent You can see that  Eq. (\ref {7}) is true in the  geometric sense, while Eq. (\ref {8}) only in the arithmetic sense. 

\noindent From (Fig.~\ref{f8}) we can conclude that Eq. (\ref{6}) is true in the geometric sense for $x,y,z,w,t>0$.

\begin{figure}[H]
\begin{center} 
\begin{tikzpicture}[>=stealth,scale=1]
\draw[fill,gray!20] (0,0)--(5,10)--(10,0)--(0,0);
\draw(0,0)--(5,10)--(10,0)--(0,0);
\draw[dashed] (0,0)--(5,2)--(10,0);
\draw[dashed] (5,2)--(5,10);
\draw[dotted,blue,thick](6,8)--(2,0)--(6,1.6)--cycle;
\draw[dotted,blue,thick](4,8)--(8,0)--(4,1.6)--cycle;
\draw[dotted,blue,thick](1,.4)--(9,.4)--(5,8.4)--cycle;
\draw[dotted,blue,thick](1,2)--(9,2)--(5,3.6)--cycle;
\draw[dotted,red,thick](3,2)--(6,3.2);\draw[dotted,red,thick](7,2)--(4,3.2);
\draw[dotted,red,thick](2,2.4)--(8,2.4);
\draw[dotted,red,thick](3,.4)--(6,6.4);\draw[dotted,red,thick](7,.4)--(4,6.4);
\draw[dotted,red,thick](5,1.2)--(5,6);
\draw[dotted,red,very thick](4,2.4)--(6,2.4)--(5,4.4)--cycle;
\draw[dotted,red,very thick](4,2.4)--(5,2.8)--(6,2.4);
\draw[dotted,red,very thick](5,2.8)--(5,4.4);
\draw[dotted,green,very thick](2,0)--(8,0);\draw[dotted,green,very thick](1,.4)--(4,1.6);
\draw[dotted,green,very thick](9,.4)--(6,1.6);\draw[dotted,green,very thick](1,2)--(4,8);
\draw[dotted,green,very thick](9,2)--(6,8);\draw[dotted,green,very thick](5,3.6)--(5,8.4);

\end{tikzpicture}
\caption{} \label{f8}
\end{center}
\end{figure}

\noindent In order to prove that Eq. ~(\ref{6}) is true in the geometric sense $\forall x,y,z,w,t \in \mathbb{R}$, we use the method  from \cite{Art}.\\
Firt we need  to make some preparations.\\

\noindent Let us note that all sixteen tetrahedrons $T_j$, $j=1,2,\ldots,16$ of equation ~(\ref{6}) arise from the intersection of four pairs of parallel planes. Indeed, the four planes $A_1, B_1, C_1,D_1$ form the tetrahedron $<t>$, while the  four planes $A_2, B_2, C_2,D_2$, parallel respectively to $A_1, B_1,C_1,D_1$, form the tetrahedron $<x,y,z,w,t>$. The faces of the remaining tetrahedrons of equation ~(\ref{7}) lie on other combinations of the planes $A_i, B_i, C_i,D_i$, where $i=1,2$.
If four non-parallel planes intersect at a single point, then some of the tetrahedra of equation ~(\ref{7}) degenerate to points.\\

\noindent Let us assume that the tetrahedra $\<x>$ are closed for $x>0$
 and open for $x<0$.\\
Each of the sixteen tetrahedra is either indivisible or consists of smaller polyhedra that are themselves indivisible. We call these smallest indivisible polyhedra (including indivisible tetrahedrons) \emph{cells} (atoms). If a polyhedron shares a face with a closed tetrahedron, we regard that face as belonging to the tetrahedron, not to the polyhedron. If a polyhedron shares a face with an open tetrahedron, we regard that face as belonging to the polyhedron. If a polyhedron lies inside  a closed tetrahedron, its remaining faces to be closed; if it lies inside  an open tetrahedron,  its remaining faces are considered open.

\noindent Let us count the number of the cells.\\ 
Three pairs of parallel planes $A_1, B_1,C_1, A_2, B_2, C_2$ form a parallelepiped (Fig.~\ref{f83}).
Through two vertices of the parallelepiped we draw the axis $l_1$, and through the six remaining vertices we draw the axes $l_2,l_3,\ldots,l_7$ , parallel to the axis $l_1$.
Then the four pairs of parallel planes form three cells $C_1, C_2, C_3$ containing the axis $l_1$, and  $6\cdot2$ cells $C_4, C_5,\ldots,C_{15}$ containing the axes $l_2,l_3,\ldots,l_7$ (only the planes $D_1, D_2$
complete the cells containing the axes $l_i$, $i\neq1$) (Fig.~\ref{f83}).\\
 
\begin{figure}[H]
\begin{center} 
\begin{tikzpicture}[>=stealth,scale=.4]
\draw[blue,fill=blue!10,dashed](1,1)--(14,7)--(24,7)--(21.2,8)--(6,1);
\draw[blue,dashed](1,1)--(18.3,9)--(24,7)--(19,7)--(6,1)--(1,1);
\draw[dashed,blue](14,7)--(19,7)--(21.2,8);
\draw[blue,fill=blue!10,dashed](13.3,19.5)--(23.3,19.5)--(21.15,20.1)--(16.75,18.4)--cycle;
\draw[blue](13.3,19.5)--(19.5,19.5)--(16.75,18.4)--cycle;
\draw[blue](1,1)--(6,1)--(19,7)--(24,7);
\draw(0,0)--(21,21)--(24,19)--(5,0);
\draw(11,21)--(20,15)--(22,3);
\draw(24,19)--(25,13)--(21,9)--(27,5);
\draw(21,9)--(16,4);
\draw[dashed](18,11)--(14,7);
\draw[dashed]((21,21)--(22,15)--(25,13);
\draw[dashed](22,15)--(18,11)--(21,9);
\draw[dashed](16.5,20)--(18.5,8);
\draw(16.5,20)--(17,17);
\draw[red,dashed,thick](21,23)--(21,2);
\draw[green,dashed,thick](17,22)--(17,5);

\draw node at(21.5,23.5){\footnotesize $l_1$};
\draw node at(17.5,22.5){\footnotesize $l_2$};
\draw node at(22,7.7){\footnotesize $C_1$};
\draw node at(22,12){\footnotesize $C_2$};
\draw node at(21.8,20){\footnotesize $C_3$};
\draw node at(14,8){\footnotesize $C_4$};
\draw node at(16.4,18){\footnotesize $C_5$};
\draw node at(12.3,5){\footnotesize $D_1$};
\draw node at(16,19.1){\footnotesize $D_2$};
\draw node at(18,7.7){\footnotesize $D_1$};
\end{tikzpicture}
\caption{Cells $C_1, C_2, C_3$ containing the axis $l_1$ and cells $C_4, C_5$ with the axis $l_2$.} \label{f83}
\end{center}
\end{figure}

\begin{thm}
      $\forall x,y,z,w,t \in \mathbb{R}$ Equation ~(\ref{6}) is true in the geometric sense. \label{t2.4}
\end{thm}
\begin{proof}
Each of the sixteenteen tetrahedra of equation ~(\ref{6}) is a uniquely defined sum of a subset of the fifteen atoms.  
We construct the matrix
\[
M \in \{0,1\}^{15 \times 16},
\]
where $M_{i,j}=1$ if the cell $C_i$ contributes to the tetrahedron $T_j$, and $0$ otherwise.  
The rank of $M$ is at most $15$. Hence, every tetrahedron is a rational linear combination of the others.  

\noindent Thus, there exist rational coefficients $A,a_i$, where $i=1,2,\ldots,15$ such that the following equation holds in the geometric sense:

\begin{align}\begin{split}
 A\<x,y,z,w,t>&= a_1\< x,y,z,0,t> + a_2\<x,y,0,w,t>\\
     &\quad + a_3\<x,0,z,w,t>+a_4\<0,y,z,w,t>  \\
     &\quad - a_5\<x,y,0,0,t> -a_6\< x,0,z,0,t> - a_7\< x,0,0,w,t > \\
     &\quad - a_8\<0,y,z,0,t>- a_9\<0,y,0,w,t>- a_{10}\<0,0,z,w,t>\\             
 &\quad + a_{11}\<x,0,0,0,t>+ a_{12}\<0,y,0,0,t>+ a_{13}\<0,0,z,0,t>\\
     &\quad +a_{14}\<0,0,0,w,t> - a_{15}\< 0,0,0,0,t >.
    \label{81}
 \end{split}   \end{align}
 
\noindent From this it follows that
\begin{align}\begin{split}
A(x+y+z+w+t)^3&= a_1( x+y+z+t)^3 + a_2(x+y+w+t)^3\\
     &\quad  + a_3(x+z+w+t)^3+a_4(y+z+w+t)^3  \\
     &\quad - a_5(x+y+t)^3 -a_6(x+z+t)^3 - a_7(x+w+t)^3\\
     &\quad - a_8(y+z+t)^3- a_9(y+w+t)^3- a_{10}(z+w+t)^3\\             
 &\quad+ a_{11}(x+t)^3+a_{12}(y+t)^3+a_{13}(z+t)^3\\
     &\quad +a_{14}(w+t)^3 - a_{15}(t)^3. 
      \label{82}       
\end{split}\end{align} 

\noindent We also know that equation ~(\ref{6}) is true in  algebraic sense. Thus,
\begin{align}\begin{split}
(x+y+z+w+t)^3&= ( x+y+z+t)^3 + (x+y+w+t)^3\\
     &\quad  + (x+z+w+t)^3+(y+z+w+t)^3  \\
     &\quad - (x+y+t)^3 -(x+z+t)^3 - (x+w+t)^3\\
     &\quad - (y+z+t)^3- (y+w+t)^3- (z+w+t)^3\\             
 &\quad+ (x+t)^3+(y+t)^3+(z+t)^3\\
     &\quad +(w+t)^3 - (t)^3.
      \label{83}        
\end{split}\end{align}

\noindent Let us compute $(\ref{81}) - a(\ref{82})$:
\begin{align}\begin{split}
(A-a_1)(x+y+z+w+t)^3&=  (a_2-a_1)(x+y+w+t)^3\\
     &\quad  + (a_3-a_1)(x+z+w+t)^3+(a_4-a_1)(y+z+w+t)^3  \\
     &\quad - (a_5-a_1)(x+y+t)^3 -(a_6-a_1)(x+z+t)^3\\
     &\quad - (a_7-a_1)(x+w+t)^3- (a_8-a_1)(y+z+t)^3\\
     &\quad - (a_9-a_1)(y+w+t)^3- (a_{10}-a_1)(z+w+t)^3\\             
     &\quad + (a_{11}-a_1)(x+t)^3+(a_{12}-a_1)(y+t)^3\\
     &\quad +(a_{13}-a_1)(z+t)^3 +(a_{14}-a_1)(w+t)^3\\
     &\quad - (a_{15}-a_1)(t)^3. 
       \label{84}       
\end{split}\end{align}

\noindent On the left-hand side of ~(\ref{84}) we have the term $(A-a)6xyz$, but it does not appear on the right-hand side. Therefore, $A=a_1$.
\noindent By a similar argument, we conclude that $A=a_2=a_3=a_4$.  

\noindent We thus obtain
\begin{align*}\begin{split}
A(x+y+z+w+t)^3&= A( x+y+z+t)^3 + A(x+y+w+t)^3\\
     &\quad  + A(x+z+w+t)^3+A(y+z+w+t)^3  \\
     &\quad - a_5(x+y+t)^3 -a_6(x+z+t)^3 - a_7(x+w+t)^3\\
     &\quad - a_8(y+z+t)^3- a_9(y+w+t)^3- a_{10}(z+w+t)^3\\             
 &\quad+ a_{11}(x+t)^3+a_{12}(y+t)^3+a_{13}(z+t)^3\\
     &\quad +a_{14}(w+t)^3 - a_{15}(t)^3.        
\end{split}\end{align*} 

\noindent But then  
\[
6Axyt = 12Axyt - 6a_5xyt,
\]
which implies $A=a_5$.  
Similarly,  $A=a_i,$, where $i=6,7,8,9,10$. 

\noindent So we obtain
\begin{align*}\begin{split}
A(x+y+z+w+t)^3&= A( x+y+z+t)^3 + A(x+y+w+t)^3\\
     &\quad  + A(x+z+w+t)^3+A(y+z+w+t)^3  \\
     &\quad - A(x+y+t)^3 -A(x+z+t)^3 - A(x+w+t)^3\\
     &\quad - A(y+z+t)^3- A(y+w+t)^3- A(z+w+t)^3\\             
     &\quad+ a_{11}(x+t)^3+a_{12}(y+t)^3+a_{13}(z+t)^3\\
     &\quad +a_{14}(w+t)^3 - a_{15}(t)^3.       
\end{split}\end{align*}  

\noindent But then  
\[
3Axyt = 9Axyt - 9Axt+ 3a_{11}xyt,
\]
which implies $A=a_{11}$.  
Similarly,  $A=a_i,$, where $i=12,13,14$.\\ 
Finally, we also obtain $A=a_{15}$. 

\noindent Therefore, equation ~(\ref{81}), true in the geometric sense, takes the form
\begin{align*}\begin{split}
A\<x,y,z,w,t>&= A\< x,y,z,0,t> + A\<x,y,0,w,t>\\
     &\quad + A\<x,0,z,w,t>+A\<0,y,z,w,t>  \\
     &\quad - A\<x,y,0,0,t> -A\< x,0,z,0,t> - A\< x,0,0,w,t > \\
     &\quad - A\<0,y,z,0,t>- A\<0,y,0,w,t>- A\<0,0,z,w,t>\\             
 &\quad + A\<x,0,0,0,t>+ A\<0,y,0,0,t>+ A\<0,0,z,0,t>\\
     &\quad +A\<0,0,0,w,t> - A\< 0,0,0,0,t >,
  \end{split}\end{align*}
for $A \neq 0$. This completes the proof. 
 
\end{proof}

\noindent The tetrahedron $\<n >$ for $n\in \mathbf{Z}$ can also be written in the following form.
\begin{align}\begin{split}
       \<n>=&\ \frac{n(n+1)(n+2)}{6}\<1> +\frac{(n-1)n(n+1)}{6}\Big(\<2>-4\<1>\Big) \\
            &\ + \frac{(n-2)(n-1)n}{2}\Big(\<3>-4\<2>+6\<1>\Big)\\
            &\ + \frac{(n-3)(n-2)(n-1)}{2}\Big(\<4>-4\<3>+6\<2>-4\<1>\Big)\\
            =&\  \frac{n(n+1)(n+2)}{6}\<1> +\frac{(n-1)n(n+1)}{6}\Big(\<2>-4\<1>\Big) \\
            &\  + \frac{(n-2)(n-1)n}{2}\Big(4\<0>-\<-1>\Big)+ \frac{(n-3)(n-2)(n-1)}{2}\big(-<0>\big).
    \label{9}
\end{split}\end{align}  

\noindent The expressions
\begin{align*}
\<a_{3,1}>=&\ \<2>-4\<1>=(4,0,-2,-3)\; \text{(Fig.~\ref{f9})} ,\\
\<a_{3,2}>=&\ \<3>-4\<2>+6\<1>=4\<0>-\<-1>=(1,-1,1,3),\\
\end{align*} 
\begin{figure}[H]
\begin{center} 
\begin{tikzpicture}[>=stealth,scale=.7]
\draw[fill,gray!20](2,0)--(4,4)--(0,4)--(2,0) ;
\draw[fill,gray!30](2,0)--(2.5,1)--(1.5,1)--cycle ;
\draw[fill,gray!30](3.25,2.5)--(2.5,4)--(4,4)--cycle ;
\draw[fill,gray!30](.75,2.5)--(1.5,4)--(0,4)--cycle ;
\draw[fill,green!20](0,4)--(4,4)--(2,5)--cycle;
\draw[fill,green!20](0,4)--(2,0)--(0,1)--cycle;
\draw[fill,green!20](4,4)--(4,1)--(2,0)--cycle;
\draw[fill,green!30](0,1)--(1.5,1)--(.75,2.5)--cycle;
\draw[fill,green!30](4,1)--(2.5,1)--(3.25,2.5)--cycle;
\draw[fill,green!30](2,5)--(2.5,4)--(1.5,4)--cycle;
\draw[green,thick](0,1)--(2,0)--(4,4)--(0,4)--(2,0)--(4,1)--(4,4)--(2,5)--(0,4)--(0,1);
\draw[dotted,green,thick](0,1)--(4,1)--(2,5)--cycle;
\draw[fill,red](0,1)circle(1pt);\draw[fill,red](2,0)circle(1pt);
\draw[fill,red](4,4)circle(1pt);\draw[fill,red](0,4)circle(1pt);
\draw[fill,red](4,1)circle(1pt);\draw[fill,red](2,5)circle(1pt);
\draw node[] at (-.3,1){\small{A}};\draw node[] at (2,-.3){\small{B}};
\draw node[] at (4.3,1){\small{C}};\draw node[] at (4.3,4){\small{D}};
\draw node[] at (2,5.3){\small{E}};\draw node[] at (-.3,4){\small{F}};
\end{tikzpicture}
\caption{The expression $\<a_{3,1}>$. The faces ABC, CDE, AEF, BDF are black.
The faces ABF, BCD, DEF, ACE are green. The octahedron is black inside.} \label{f9}
\end{center}
\end{figure}

\noindent not belongs to  $\mathbb{R}_{3}$ but they are useful for a demonstration of the geometric interpretation of the tetrahedra (Fig.~\ref{f1}).\\
We can see that the first coordinates of the expressions $\<a_{3,0}>=\<1>=(1,1,1,1) $,  $\<a_{3,1}> $ and $\<a_{3,2}>$,  are Eulerian numbers: $A_{3,0}=1, \; A_{3,1}=4, \; A_{3,2}=1$ \cite{GKP}.\\
But when we add the expression\\
 $\<a_{3,3}>=\<4>-4\<3>+6\<2>-4\<1>=-\<0>=(0,0,0,-1) $,
  then all the expressions $\<a_{3,i}>, \; i=0,1,2,3$ can be treated as a generalization of Eulerian numbers.
 Let us note that \\
 \[ \sum_{i=0}^3\<a_{3,i}>=(6,0,0,0). \] 
\noindent (Fig.~\ref{f10})

\begin{figure}[H]
\begin{center} 
\begin{tikzpicture}[>=stealth,scale=.5]
\draw[fill,gray!30] (-5,4)--(-1,4)--(-3,8)--cycle;
\draw (-5,4)--(-1,4)--(-3,8)--cycle;
\draw[dashed](-5,4)--(-3,5)--(-1,4);\draw[dashed](-3,5)--(-3,8);
\node at(-.5,4){+};

\draw[fill,gray!20](2,0)--(4,4)--(0,4)--(2,0) ;
\draw[fill,green!20](0,4)--(4,4)--(2,5)--cycle;
\draw[fill,green!20](0,4)--(2,0)--(0,1)--cycle;
\draw[fill,green!20](4,4)--(4,1)--(2,0)--cycle;
\draw[green,thick](0,1)--(2,0)--(4,4)--(0,4)--(2,0)--(4,1)--(4,4)--(2,5)--(0,4)--(0,1);
\draw[dotted,green,thick](0,1)--(4,1)--(2,5)--cycle;
\draw[fill,red](0,1)circle(1.5pt);\draw[fill,red](2,0)circle(1.5pt);
\draw[fill,red](4,4)circle(1.5pt);\draw[fill,red](0,4)circle(1.5pt);
\draw[fill,red](4,1)circle(1.5pt);\draw[fill,red](2,5)circle(1.5pt);
\node at(4.5,1){+};

\draw[fill,gray!30] (5,1)--(9,1)--(7,-3)--cycle;
\draw[green,thick] (5,1)--(9,1)--(7,-3)--(5,1)--(7,0)--(9,1);
\draw[green,thick](7,0)--(7,-3);
\draw[fill](5,1)circle(1.5pt);\draw[fill](9,1)circle(1.5pt);
\draw[fill](7,-3)circle(1.5pt);\draw[fill](7,0)circle(1.5pt);
\node at(9.5,-3){+};

\draw[fill,red](10,-3)circle(1.5pt);
\node at(10.5,2.5){=};

\draw[fill,gray!20](11,4)--(13,0)--(15,4)--(13,8)--cycle;
\draw(11,4)--(13,8)--(15,4);\draw[dashed](13,8)--(13,5);
\draw[fill,green!20](11,4)--(13,0)--(15,4)--(15,1)--(13,-3)--(11,1)--cycle;
\draw[green,thick](11,4)--(13,0)--(13,-3)--(11,1)--cycle;
\draw[green,thick](13,0)--(15,4)--(15,1)--(13,-3)--cycle;
\draw[green,dashed,thick](11,1)--(13,5)--(15,1);
\draw[fill,green](11,4)circle(1.5pt);\draw[fill,green](13,5)circle(1.5pt);
\draw[fill,green](15,4)circle(1.5pt);\draw[fill,green](11,1)circle(1.5pt);
\draw[fill,green](13,0)circle(1.5pt);\draw[fill,green](15,1)circle(1.5pt);
\draw[fill,green](13,-3)circle(1.5pt);
\draw[fill](13,8)circle(1.5pt);

\end{tikzpicture}
\caption{$\<a_{3,2}>$ has green faces and is black inside. 
The top three faces of solid (6,0,0,0) are black and the three bottom faces are green. The solid is black inside.} \label{f10}
\end{center}
\end{figure}
  
\noindent Eq. (\ref{5}) is special case of the following equation
\begin{align}\begin{split}
\forall x,a,b,c,d \in \mathbb{R}\qquad & \\ 
       \<x>=&\ \frac{(x-b)(x-c)(x-d)}{(a-b)(a-c)(a-d)}\<a>-\frac{(x-a)(x-c)(x-d)}{(b-a)(b-c)(b-d)}\<b>\\
        &\ +\frac{(x-a)(x-b)(x-d)}{(c-a)(c-b)(c-d)}\<c> + \frac{(x-a)(x-b)(x-c)}{(d-a)(d-b)(d-c)}\<d> .
\label{10}
\end{split} \end{align} 
\noindent where coefficients are the Lagrange interpolating polynomials.\\
Eq. (\ref{10}) is true in the arithmetic sense but is not true in the geometric sense for any $x,\ a,\ b,\ c \in \mathbb{R} $.\\
You can investigate for which values $x,a,b,c,d $ Eq. (\ref{10})  is true in the geometric sense.
For example the following  equation is probably true in the geometric sense for all $n,k \in \mathbb{Z}$.
\begin{align*}
       \<n>=&\ \frac{(n-k+1)(n-k+2)(n-k+3)}{6}\<k>-\frac{(n-k)(n-k+2)(n-k+3)}{2}\<k-1>\\
        &\ +\frac{(n-k)(n-k+1)(n-k+3)}{2}\<k-2> - \frac{(n-k)(n-k+1)(n-k+2)}{6}\<k-3> .
 \end{align*} 
In general, the criterion of whether a given equation is true in the geometric sense is unknown.\\

\section{Multiplication in the set $\mathbb{R}_3$.}

\noindent Let us fix a point in $\mathbb{R}_3$ and denote it as $\<0,0,0,0>$. Then, using Eq. (\ref{6}), every element of the set $\mathbb{R}_3$ can be obtained from the distinguished element $\<0,0,0,0>$ as an element $\<x,y,z,w,0>$. We will mark it $\<x,y,z,w>$.\\
For the element $\<t>=\<t_x,t_y,t_z,t_w,0>=\<t_x,t_y,t_z,t_w> $
 Eq. (\ref{6}) takes the following form
 
 \begin{align}\begin{split}
  \forall x,y,z,w,t \in \mathbb{R}  \\
   \<x,y,z,w,t>&= \< x+t_x,y+t_y,z+t_z,w+t_w>\\       &=\< x+t_x,y+t_y,z+t_z,t_w> + 
         \<x+t_x,y+t_y,t_z,w+t_w>\\ 
& \quad + \<x+t_x,t_y,z+t_z,w+t_w>  +
          \<t_x,y+t_y,z+t_z,w+t_w>\\
&\quad - \<x+t_x,y+t_y,t_z,t_w>-\<x+t_x,t_y,z+t_z,t_w>           
       -  \< x+t_x,t_y,t_z,w+t_w > \\
&\quad - \<t_x,y+t_y,z+t_z,t_w>-\<t_x,y+t_y,z+t_z,t_w>           
       -  \< t_x,t_y,z+t_z,w+t_w > \\
&\quad + \<x+t_x,t_y,t_z,t_w>+\<t_x,y+t_y,t_z,t_w>           
       +  \< t_x,t_y,z+t_z,t_w > \\       
&\quad + \< t_x,t_y,t_z,w+t_w >- \< t_x,t_y,t_z,t_w >.    \label{101}
 \end{split}   \end{align}

\noindent It is easy to see that the elements $\<x,y,z,w>$ can be interpreted as points of the space $\mathbb{R}^4$.\\
So all triangles from   $\mathbb{R}_3$ satisfying the condition $x+y+z+w=t$ for a fixed $t$, correspond to a hyperplane in $\mathbb{R}^4$.\\
Let us also observe that the elements of Eq. (\ref{101}) form the vertices of a  four-dimensional rectangular cuboid  in $\mathbb{R}^4$.\\

\subsection{Multiplication in the set $\mathbb{R}_3$}

\noindent We know that the set $\mathbb{R}_3$ is closed under multiplication. So what does multiplication of elements $\<x,y,z,w> \in \mathbb{R}_3$ look like? Let us assume that all components of $\<x,y,z,w> $ are equivalent,  multiplication is associative and commutative, and multiplication is distributive over addition with respect to the components.  Then the following equation must hold.  
\begin{align}
     \Big(\<1,0,0,0>^2\Big)\<0,1,0,0>=\<1,0,0,0>\Big(\<1,0,0,0>\<0,1,0,0>\Big).
     										\label{102}
 \end{align}   

\noindent Let us suppose that \\
$\<1,0,0,0>^2= \<1-3a,a,a,a>, \qquad
   \<1,0,0,0>\<0,1,0,0>=\<c,c,\frac{1}{2}-c,\frac{1}{2}-c>$.\\
Then
\begin{align*}\begin{split}
 \Big(\<1,0,0,0>^2\Big)& \<0,1,0,0>= \<1-3a,a,a,a>\<0,1,0,0>  \\
        = &\ \<1-3a,0,0,0>\<0,1,0,0> +\<0,a,0,0>\<0,1,0,0>\\
          &\ +\<0,0,a,0>\<0,1,0,0>+\<0,0,0,a>\<0,1,0,0>\\
        = &\ \<c-3ac,c-3ac,(1-3a)\Big(\frac{1}{2}-c\Big),
                     (1-3a)\Big(\frac{1}{2}-c\Big)>\\ 
          &\ +\<a^2,a-3a^2,a^2,a^2>\\
          &\ +\<\frac{a}{2}-ac,ac,ac,\frac{a}{2}-ac>\\
          &\ +\<\frac{a}{2}-ac,ac,\frac{a}{2}-ac,ac>\\
        = &\ \langle a+c-5ac+a^2,a+c-ac-3a^2,\\
    & \qquad\qquad  \frac{1}{2}-a-c+3ac+a^2,\frac{1}{2}-a-c+3ac+a^2\rangle             
 \end{split}   \end{align*}
\noindent and
\begin{align*}\begin{split}
\<1,0,0,0>& \Big(\<1,0,0,0>\<0,1,0,0>\Big) = 
                \<1,0,0,0>\<c,c,\frac{1}{2}-c,\frac{1}{2}-c>  \\
        = &\ \<1,0,0,0>\<c,0,0,0> +\<1,0,0,0>\<0,c,0,0>\\
          &\ +\<1,0,0,0>\<0,0,\frac{1}{2}-c,0>+
          	     \<1,0,0,0>\<0,0,0,\frac{1}{2}-c>\\
        = &\ \<(1-3a)c,ac,ac,ac>\\ 
          &\ +\<c^2,c^2,\frac{c}{2}-c^2,\frac{c}{2}-c^2>\\
          &\ +\langle\frac{c}{2}-c^2,\Big(\frac{1}{2}-c\Big)^2,
              \frac{c}{2}-c^2,\Big(\frac{1}{2}-c\Big)^2\rangle\\
          &\ +\langle\frac{c}{2}-c^2,\Big(\frac{1}{2}-c\Big)^2, 
              \Big(\frac{1}{2}-c\Big)^2,\frac{c}{2}-c^2\rangle\\                                  
        = &\ \langle 2c-3ac-c^2,\frac{1}{2}+ac-2c+3c^2,
               \frac{1}{4}+ac-c^2,\frac{1}{4}+ac-c^2\rangle. 
 \end{split}   \end{align*}
 
 \noindent According to Eq. (\ref{102}) we get
 \begin{multline}
 \< a+c-5ac+a^2,a+c-ac-3a^2,
   \frac{1}{2}-a-c+3ac+a^2,\frac{1}{2}-a-c+3ac+a^2>\\
 = \<2c-3ac-c^2,\frac{1}{2}+ac-2c+3c^2,
               \frac{1}{4}+ac-c^2,\frac{1}{4}+ac-c^2>.
     										\label{103}
 \end{multline}  
 \noindent If we compare the first components of the elements in the Eq. (\ref{103}), we get two cases\\
 $a=c=\frac{1}{3} \quad$ or $\quad c=a+1$.\\
  The first case appears to be of little significance, therefore we shall consider the second case.\\
 Let us compare the second components of the elements in the Eq. (\ref{103}) and put $\quad c=a+1$.\\
 We get $\quad a=-\frac{1}{4}, \qquad c= \frac{3}{4}$.\\
 
\noindent It is easy to check that  the obtained values satisfy  Eq.(\ref{103}) for the third and fourth components. So
\begin{align*}
     \<1,0,0,0>^2=\Big\langle\frac{7}{4},-\frac{1}{4},-\frac{1}{4},       
          -\frac{1}{4} \Big\rangle, \qquad	
    \<1,0,0,0>\<0,1,0,0>=\Big\langle\frac{3}{4},\frac{3}{4},-\frac{1}{4},-\frac{1}{4}\Big\rangle.	 		
 \end{align*}   
We assumed that all components of $\<x,y,z,v> $ are equivalent therefore
\begin{align*}
    \<0,1,0,0>^2 & =\Big\langle-\frac{1}{4},\frac{7}{4},-\frac{1}{4},
    -\frac{1}{4}\Big\rangle, & \quad	
        \<1,0,0,0>\<0,0,1,0> & =\Big\langle\frac{3}{4},-\frac{1}{4},
           \frac{3}{4},-\frac{1}{4}\Big\rangle, \\[0.2cm]  
     \<0,0,1,0>^2 & =\Big\langle-\frac{1}{4},-\frac{1}{4},\frac{7}{4},
       -\frac{1}{4}\Big\rangle, & \quad
          \<1,0,0,0>\<0,0,0,1> & =\Big\langle\frac{3}{4},-\frac{1}{4},- 
              \frac{1}{4},\frac{3}{4}\Big\rangle,\\[0.2cm]
     \<0,0,0,1>^2 & =\Big\langle-\frac{1}{4},-\frac{1}{4},-\frac{1}{4},
        \frac{7}{4}\Big\rangle, & \quad 
            \<0,1,0,0>\<0,0,1,0> & =\Big\langle-\frac{1}{4},\frac{3}{4}, 
              \frac{3}{4},-\frac{1}{4}\Big\rangle,
\end{align*}
\[    
    \<0,1,0,0>\<0,0,0,1>  =\Big\langle-\frac{1}{4},-\frac{1}{4}, 
              -\frac{1}{4},\frac{3}{4}\Big\rangle,
     \;   \<0,0,1,0>\<0,0,0,1> =\Big\langle-\frac{1}{4},-\frac{1}{4}, 
              \frac{3}{4},\frac{3}{4}\Big\rangle,  
 \]  
\\ 
   
\noindent Let us find a general formula for multiplication for arbitrary $x,y,z,v,a,b,c,d \in \mathbb{R}$.
\begin{multline}
\<x,y,z,v>\<a,b,c,d> \\
   =\Big\langle\frac{7xa}{4},-\frac{xa}{4},-\frac{xa}{4},-\frac{xa}{4}    
          \Big\rangle
  +\Big\langle\frac{3xb}{4},\frac{3xb}{4},-\frac{xb}{4},-\frac{xb}{4}
               \Big\rangle     
  +\Big\langle\frac{3xc}{4},-\frac{xc}{4},\frac{3c}{4},-\frac{xc}{4}
                     \Big\rangle\\[0.1cm]
   +\Big\langle\frac{3xd}{4},-\frac{xd}{4},-\frac{xd}{4},\frac{3xd}{4}
          \Big\rangle
  +\Big\langle\frac{3ya}{4},\frac{3ya}{4},-\frac{ya}{4},-\frac{ya}{4}
               \Big\rangle
  +\Big\langle-\frac{yb}{4},\frac{7yb}{4},-\frac{yb}{4},-\frac{yb}{4}
                    \Big\rangle\\[0.1cm]
  +\Big\langle-\frac{yc}{4},\frac{3yc}{4},\frac{3yc}{4},-\frac{yc}{4} 
          \Big\rangle
  +\Big\langle-\frac{yd}{4},\frac{3yd}{4},-\frac{yd}{4},\frac{3yd}{4} 
               \Big\rangle        
  +\Big\langle\frac{3za}{4},-\frac{za}{4},\frac{3za}{4},-\frac{za}{4}
                    \Big\rangle\\[0.1cm]
  +\Big\langle-\frac{zb}{4},\frac{3zb}{4},\frac{3zb}{4},-\frac{zb}{4}
          \Big\rangle
  +\Big\langle-\frac{zc}{4},-\frac{zc}{4},\frac{7zc}{4},-\frac{zc}{4}
                \Big\rangle
  +\Big\langle-\frac{zd}{4},-\frac{zd}{4},\frac{3zd}{4},\frac{3zd}{4}
                    \Big\rangle\\[0.1cm]
  +\Big\langle\frac{wa}{4},-\frac{wa}{4},-\frac{wa}{4},\frac{3wa}{4}
          \Big\rangle
  +\Big\langle-\frac{wb}{4},\frac{3wb}{4},-\frac{wb}{4},\frac{wb}{4}
                \Big\rangle
  +\Big\langle-\frac{wc}{4},-\frac{wc}{4},\frac{3wc}{4},\frac{3wc}{4}
                    \Big\rangle\\[0.1cm]    
  +\Big\langle-\frac{wd}{4},-\frac{wd}{4},-\frac{wd}{4},\frac{7wd}{4}
                    \Big\rangle\\[0.1cm]
   =  \Big\langle \frac{(7x-y-z-w)(a+b+c+d)+(x+y+z+w)(7a-b-c-d)}{8},\\
             \frac{(-x+7y-z-w)(a+b+c+d)+(x+y+z+w)(-a+7b-c-d)}{8},\\
          \frac{(-x-y+7z-w)(a+b+c+d)+(x+y+z+w)(-a-b+7c-d)}{8}\\
          \frac{(-x-y-z+7w)(a+b+c+d)+(x+y+z+w)(-a-b-c+7d)}{8}\Big\rangle            .  \label{104}
\end{multline}  
 
\noindent Now we can verify that the associative property of multiplication holds for any elements of the set $\mathbb{R}_3$.

\begin{thm}
\begin{align*} 
 \forall x,y,z,w,a,b,c,d,A,B,C,D \in \mathbb{R} \qquad\qquad\\
 \Big(\<x,y,z,w>\<a,b,c,d>\Big)\<A,B,C,D>&\ = \<x,y,z,w>\Big(\<a,b,c,d>\<A,B,C,D>\Big)
 \end{align*}
 \end{thm}
 \begin{proof}
Becouse  all components of $\<x,y,z,w> $ are equivalent it is enough to verify associativity on first components of $\<x,y,z,w>$.\\
Let us denote
\[\<x,y,z,w>\<a,b,c,d>=\<f,g,h,k>,\]
where
\[
f=  \frac{(7x-y-z-w)(a+b+c+d)+(x+y+z+w)(7a-b-c-d)}{8},\]
\begin{multline*}
-g-h-k \\
     = - \frac{(-3x+5y+5z+5w)(a+b+c+d)+(x+y+z+w)(-3a+5b+5c+5d)}{8}\\ 
\end{multline*}      
Then
\begin{multline*}
\Big(\<x,y,z,w>\<a,b,c,d>\Big)\<A,B,C,D>  = \<f,g,h,k>\<A,B,C,D>\\
  = \Big\langle \frac{(7f-g-h-k)(A+B+C+D)+(f+g+h+k)(5A-B-C-D)}{8},\dots
      \Big\rangle\\
\end{multline*}
We have
\[
7f-g-h-k  = (7x-y-z-w)(a+b+c+d)+(x+y+z+w)(6a-2b-2c-2d)  
\]
and
\[f+g+h+k= (x+y+z+w)(a+b+c+d)\]
The numerator of the first component of the expression \\
$\Big(\<x,y,z,w>\<a,b,c,d>\Big)\<A,B,C,D> $ will take the form
\begin{multline*}
(7f-g-h-k)(A+B+C+D)+(f+g+h+k)(7A-B-C-D)\\
  =\big[(7x-y-z-w)(a+b+c+d) +(x+y+z+w)(6a-2b)\big](A+B+C+D) \\
     +(x+y+z+v)(a+b+c+d)(7A-B-C-D)
\end{multline*}
If we substitute the expression $(7x-y-z-w)$ with $(7A-B-C-D)$, 
$\<x+y+z+w> $ with $\<A+B+C+D>$ and vice versa, the above equation remains unchanged. This shows that multiplication is associtive.\\
 \end{proof}

\noindent From Eq. (\ref{104}) we get
 \begin{multline*}
 \<x,y,z,w>^2 = \\
  \Big\langle \frac{(x+y+z+w)(7x-y-z-w)}{4},
  \frac{(x+y+z+w)(-x+7y-z-w)}{4},\\
      \frac{(x+y+z+w)(-x-y+5z+w)}{4},
       \frac{(x+y+z+w)(-x-y-z-7w)}{4}\Big\rangle.
\end{multline*}
 
\noindent From the equation 
\[ \Big\langle \frac{1}{4},\frac{1}{4},\frac{1}{4},\frac{1}{4}\Big\rangle                  
\<x,y,z,w> = \<x,y,z,w>\]
it follows that the element $\big\langle \frac{1}{4},\frac{1}{4},\frac{1}{4},\frac{1}{4}\big\rangle $  is the multiplicative identity with respect to the multiplication defined by Eq. (\ref{104}).\\

\noindent As shown in the next equation
\[\<x,y,z,-x-y-z>\<a,b,c,-a-b-c>=\<0,0,0,0>,\]
the product of two points is the  point $\<0,0,0,0>$.\\

\noindent \subsection{Roots in the set $\mathbb{R}_3$}

\noindent It is straightforward to derive formulas for the roots of elements 
 $\<x,y,z+w>$, where $x+y+z+w \neq 0$.\\
 \begin{multline*}
\sqrt{\<x,y,z,w>} \\ = 
  \Big\langle\frac{5x+y+z+w}{8\sqrt{x+y+z+w}},
  \frac{x+4y+z+w}{8\sqrt{x+y+z+w}},
  \frac{x+y+4z+w}{8\sqrt{x+y+z+w}}
  \frac{x+y+z+4w}{8\sqrt{x+y+z+w}}\Big\rangle,
\end{multline*}
\[  
\sqrt{\<x,y,z,1-x-y-z>} = 
  \Big\langle\frac{4x+1}{8},\frac{4y+1}{8},\frac{4z+1}{8},
  \frac{5-4x-4y-4z}{8}\Big\rangle.
\]
\\

\subsection{Division in the set $\mathbb{R}_3$}

For all 
$  x,y,z,w,a,b,c,d \in \mathbb{R}$, such that  $a+b+c+d\neq 0 $
 there exist $D,E,F,G \in \mathbb{R} $ such that
\begin{align}
\<x,y,z,w>= \<D,E,F,G>\<a,b,c,d>, \label{106}
\end{align}
 where\\
 \begin{align*}
  D & = \frac{x(a+5b+5c+5d)+ 
           (y+z+w)(-3a+b+c+d)}{8(a+b+c+d)^2}\\[0.1cm]
    & = \frac{(5x-3y-3z-3w)(a+b+c+d)+
           (x+y+z+w)(-3a+5b+5c+5d)}{8(a+b+c+d)^2},    
\end{align*}
\begin{align*}
  E & = \frac{y(5a+b+5c+5d)+ 
           (x+z+w)(a-3b+c+d)}{8(a+b+c+d)^2},     
\end{align*}
\begin{align*}
  F & = \frac{z(5a+5b+c+5d)+ 
           (x+y+w)(a+b-3c+d)}{8(a+b+c+d)^2},           
\end{align*}
\begin{align*}
  G & = \frac{w(5a+5b+5c+d)+ 
           (x+y+z)(a+b+c-3d)}{8(a+b+c+d)^2}.            
\end{align*}
So if an element $\<a,b,c,d>$ is not a point, then it divides 
any  tetrahedron $\<x,y,z,w>$.\\

\section{The set $\mathbb{R}_k$\ .}  
  
Let us take the ring $\mathbb{R}^{k+1}
= \{ (x_1,x_2,\ldots,x_{k+1}) ; \; x_1,x_2,\ldots,x_{k+1}\in\mathbb{R} \}$ with  additon and multiplication\\
\begin{align}
(x_1,x_2,\ldots,x_{k+1})+(y_1,y_2,\ldots,y_{k+1})=&\ (x_1+y_1,x_2+y_2,\ldots,x_{k+1}+y_{k+1}), \label{11} \\
(x_1,x_2,\ldots,x_{k+1})\cdot(y_1,y_2,\ldots,y_{k+1}) =&\ (x_1\cdot y_1,x_2\cdot y_2,\ldots,x_{k+1}\cdot y_{k+1}). \label{12}
\end{align}
Let us consider the subset of the ring $\mathbb{R}^{k+1}$,\\
 the set $\mathbb{R}_k= \{\pm  \<x >=\pm (x^k,x^{k-1},\ldots,x,1); x\in\mathbb{R} \}$.\\  
It is closed under multiplication (\ref{12}) but not under addition (\ref{11}).\\
The set $\mathbb{R}_k$ is closed under    the following kind of addition 
(It should be noted that the summands in below Eq.(\ref{13}) are not  the coordinates of the ring $\mathbb{R}^{k+1})$.
 
\begin{align}\begin{split}
\forall x_1,x_2,\ldots,x_{k+1} \in \mathbb{R}\qquad & \\ 
       \<x_1+x_2+\ldots+x_{k+1}+t> = &\ \sum_{j=1}^{k+1} \< \sum_{\substack{m=1\\ m\neq j}}^{k+1} x_m+t >
      - \sum_{\substack{i,j=1\\ i\neq j}}^{k+1}	\< \sum_{\substack{m=1\\ n\neq i,m\neq j}}^{k+1} x_m+t >\\   
      &\ + \cdots +(-1)^k \sum_{\substack{m,n=1\\ m\neq n}}^{k+1}\<x_m+x_n+t> \\
      &\  	+(-1)^{k+1}\sum_{m=1}^{k+1}\<x_m+t> +(-1)^k \<t> \label{13}
\end{split}\end{align}
bacause the equations 
\begin{align*}
\forall x_1,x_2,\ldots,x_{k+1} \in \mathbb{R}\qquad \qquad &\ \forall p=k,k-1,\ldots,1,0  \\ 
       (x_1+x_2+\ldots+x_{k+1}+t)^p = &\ \sum_{j=1}^{k+1} \Big( \sum_{\substack{m=1\\ m\neq j}}^{k+1} x_m+t\Big)^p
      - \sum_{\substack{i,j=1\\ i\neq j}}^{k+1}	\Big( \sum_{\substack{m=1\\ n\neq i,m\neq j}}^{k+1} x_m+t \Big)^p\\   
      &\ + \cdots +(-1)^k \sum_{\substack{m,n=1\\ m\neq n}}^{k+1}(x_m+x_n+t)^p \\
      &\  	+(-1)^{k+1}\sum_{m=1}^{k+1}(x_m+t)^p +(-1)^k (t)^p 
\end{align*}
are true.
If we multiply the Eq. (\ref{13}) by $-1$, we get the definition of addition for the elements $-\<x>$ of the set $\mathbb{R}_k$.\\
It is easy to check that the following equation holds for each $x\in \mathbb{R}$.
\begin{align}\begin{split}
       \<x>=&\ \frac{(x-k+2)(x-k+3)\cdot\ldots\cdot (x-1)x(x+1)}{k!}\<k-1> \\
           &\ -\frac{(x-k+1)(x-k+3)\cdot\ldots\cdot (x-1)x(x+1)}{(k-1)!}\<k-2>\\
           &\ +\frac{(x-k+1)(x-k+2)(x-k+4)\cdot\ldots\cdot (x-1)x(x+1)}{2(k-2)!}\<k-3>\\
           &\ -\cdots+(-1)^k\frac{(x-k+1)(x-k+2)\cdot\ldots\cdot (x-2)x(x+1)}{2(k-2)!}\<1>\\
           &\  +(-1)^{k+1}\frac{(x-k+1)(x-k+2)\cdot\ldots\cdot (x-2)(x-1)(x+1)}{(k-1)!}\<0>\\
           &\ +(-1)^k\frac{(x-k+1)(x-k+2)\cdot\ldots\cdot (x-2)(x-1)x}{k!}\<-1>.\\
    \label{14}
\end{split}\end{align} 
To prove the truth of the Eq. (\ref{14}), it is enough to check it for the values of $k-1,k-2,\ldots,1,0,-1$.\\

\section{Geometric interpretation of the set $\mathbb{R}_k$\ .}

\noindent  Let us set any closed k dimensional  simplex   (the simplex does not have to be regular) and let us denote it by the symbol $\< 1>$. Each simplex obtained from the fixed simplex $\<1>$ by any translation will be denoted by $\<1>$ too.\\
The simplex obtained from the simplex $\<1>$ by any homothety of ratio $x>0$ \cite{Wol} will be denoted by $\<x>$.\\ 
Each point of the k dimensional space will be denoted by $\<0>$.\\
The symbol  $-\< x>$ denotes the simplex, which if is a cover of  the simplex  $\< x>$ in topological sense  \cite{Wol} and vice versa  the simplex  $\< x>$ is a cover of  the simplex  $-\< x>$ then both give an empty set  denoted by an element $(0,0,\ldots,0) \in \mathbb{R}^{k+1}$.\\
We will write easier that if the simplex  $-\< x>$ lies on the the simplex  $\< x>$ then both give an empty set.
We will  call the simplex   $-\< x>$ sometimes antysimplex .\\
An element $\<-x>$ for each $x \in \mathbb{R},\: x>0$  will be interpreted as a open simplex symmetrical to the  simplex  $\<x >$ in relation to any edge. \\
At the end an element $-\<-x>$ for each $x \in \mathbb{R},\: x>0$  will be interpreted as an open simplex symmetrical to the  simplex  $-\<x >$ in relation to any edge.\\

\noindent  Geometric interpretation of the adding (\ref{13}) is similar to adding in the set $\mathbb{R}_3$.  
\\ 
 \noindent Let us set the common order $d_1,d_2,\ldots,d_{k+1}$ of $k+1$   faces of  the each simplex $\<t>$. 
 Let us fix the ordered successive components  $ x_1,x_2,\ldots,x_{k+1}  $ of the sum \mbox{$ \< x_1+x_2+\ldots+x_{k+1}+t > $.} This components extend or shorten the simplex $ \< t> $  towards the faces \mbox{$d_1,d_2,\ldots,d_{k+1}$ }  depending on whether the numbers $x_1,x_2,\ldots,x_{k+1} $     are positive or negative (Fig.~\ref{f2}).\\
 
 \noindent As in the set $\mathbb{R}_3$, we can write definitions.
\begin{definition}
The equation $\< x> = \sum_j\alpha_j\<x_j>$, where $\alpha_j \in \{-1,1\}$ is true in the geometric sense if we can build the simplex $\< x>$ from the elements $\alpha_j\< x_j>$. \label{d3.1}
\end{definition}

\begin{definition}
The equation $\< x> = \sum_j\alpha_j\<x_j>$, where $\alpha_j \in \{-1,1\}$ is true in the arithmetic sense if the equations  $ x^i  = \sum_j\alpha_j x^i_j$, where $\alpha_j \in \{-1,1\}$,  $i=0,1,2,\ldots,k$ \label{d3.2}
hold.
\end{definition}

\noindent In order to prove that Eq. ~(\ref{13}) is true in the geometric sense, we use the method  from the proof of the theorem ~\ref{t2.4}.\\

\noindent Let us note that all $2^{k+1}$ simplices $T_j$, $j=1,2,\ldots,2^{k+1}$ of equation ~(\ref{13}) arise from the intersection of  $k+1$ pairs of parallel hiperplanes. Indeed, the $k+1$ planes $A_i$, $i=1,2,\ldots,k+1$  form the simplex $<t>$, while the  $k+1$ hiperplanes $B_i$, $i=1,2,\ldots,k+1$, parallel respectively to $A_i$, $i=1,2,\ldots,k+1$, form the  simplex $\<x_1+x_2+\ldots+x_{k+1}+t>$. The facets of the remaining simplices of equation (~\ref{13}) lie on other combinations of the planes $A_i, B_i$, where $i=1,2,\ldots,k+1$.
If $k+1$ non-parallel hiperplanes intersect at a single point, then some of the simpices of equation ~(\ref{13}) degenerate to points.\\

\noindent Let us assume that the simplices $\<x>$ are closed for $x>0$
 and open for $x<0$.\\
Each of the $2^{k+1}$ simplces is either indivisible or consists of smaller polytopes that are themselves indivisible. We call these smallest indivisible polytopes (including indivisible simpices) \emph{cells} (atoms). If a polytope shares a facet with a closed simplex, we regard that facet as belonging to the simplex, not to the polytope. If a polytope shares a facet with an open tetrahedron, we regard that facet as belonging to the polytope. If a polytope lies inside  a closed simplex, its remaining facets to be closed; if it lies inside  an open simplex,  its remaining facets are considered open.

\noindent Let us count the number of the cells.\\ 
$k$ pairs of parallel planes $A_i, B_i$, $i=1,2,\ldots,k$ form a parallelotope.
Through two vertices of the paralleotope we draw the axis $l_1$, and through the $2^k-2$ remaining vertices we draw the axes $l_2,l_3,\ldots,l_{2^k-1}$ , parallel to the axis $l_1$.
Then the all pairs of parallel hiperplanes form three cells $C_1, C_2, C_3$ containing the axis $l_1$, and  $(2^k-2)\cdot2$ cells $C_4, C_5,\ldots,C_{2^{k+1}-1}$ containing the axes $l_2,l_3,\ldots,l_{2^k-1}$ (only the hiperplanes $D_1, D_2$
complete the polytopes containing the axes $l_i$, $i\neq1$).\\

\begin{thm}
 $\forall x_1,x_2,\ldots,x_{k+1} \in \mathbb{R}$
    Equation (\ref{13}) is true in the geometric sense. \label{t5.3}
\end{thm}
\begin{proof}[sketch of the proof]

Each of the $2^{k+1}$ simplices of equation (~\ref{13}) is a uniquely defined sum of a subset of the $2^{k+1}-1$ atoms.  
We construct the matrix
\[
M \in \{0,1\}^{(2^{k+1}-1) \times (2^{k+1})},
\]
where $M_{i,j}=1$ if the cell $C_i$ contributes to the simpices $T_j$, and $0$ otherwise.  
The rank of $M$ is at most $2^{k+1}-1$. Hence, every simplex is a rational linear combination of the others.  

\noindent Thus, there exist rational coefficients $A,a_{ij}$  such that the following equation holds in the geometric sense:

\begin{align*}\begin{split}
   A\<x_1+x_2+\ldots+x_{k+1}+t> = &\ \sum_{j=1}^{k+1}a_{1j} \< \sum_{\substack{m=1\\ m\neq j}}^{k+1} x_m+t >
      - \sum_{\substack{i,j=1\\ i\neq j}}^{k+1}a_{2j}	\< \sum_{\substack{m=1\\ n\neq i,m\neq j}}^{k+1} x_m+t >\\   
      &\ + \cdots +(-1)^k a_{(k-1)j}\sum_{\substack{m,n=1\\ m\neq n}}^{k+1}\<x_m+x_n+t> \\
      &\  	+(-1)^{k+1}a_{kj}\sum_{m=1}^{k+1}\<x_m+t> +(-1)^ka_{(k+1)j} \<t> 
   \end{split}   \end{align*}
 By applying the same operations as in the theorem ~\ref{t2.4}, we conclude that all the coefficients $A,a_{ij}$
 are equal. This completes the proof.  
 \end{proof}

\noindent The following equation is  true in the arithmetic sense for each $n\in \mathbb{Z}$.
\begin{align}\begin{split}
       \<n>=&\ \frac{n(n+1)\cdot\ldots\cdot (n+k-1)}{k!}\<1> \\
            &\ +\frac{(n-1)n\cdot\ldots\cdot (n+k-2)}{k!}\Bigg(\<2>-{k+1\choose 1}\<1>\Bigg) \\
            &\ +\frac{(n-2)(n-1)\cdot\ldots\cdot (n+k-3)}{k!}
                  \Bigg(\<3>-{k+1\choose 1}\<2>+{k+1\choose 2}\<1>\Bigg) \\
            &\ +\cdots +\frac{(n-k+1)(n-k+2)\cdot\ldots\cdot n}{k!}\\
            &\ \cdot \Bigg(\<k>-{k+1\choose 1}\<k-1>+{k+1\choose 2}\<k-2>- 
              \cdots+(-1)^{k+1}{k+1\choose k-1}\<1>  \Bigg) \\
           &\ +\frac{(n-k)(n-k+1)\cdot\ldots\cdot (n-1)}{k!}\\
            &\ \cdot \Bigg(\<k+1>-{k+1\choose 1}\<k>+{k+1\choose 2}\<k-1>- 
              \cdots+(-1)^k{k+1\choose k}\<1>  \Bigg). \\         
    \label{15}
\end{split}\end{align} 
Probably Eq. (\ref{15}) is true in the geometric sense.\\
\noindent The expressions
\begin{align*}
\<a_{k,0}>=&\ \<1>,\\
\<a_{k,1}>=&\ \<2>-{k+1\choose 1}\<1>,\\
\<a_{k,2}>=&\ \<3>-{k+1\choose 1}\<2>+{k+1\choose 2}\<1>, \\
...........& ............................................\\
\<a_{k,k-1}> =&\ \<k>-{k+1\choose 1}\<k-1>+{k+1\choose 2}\<k-2>- \cdots+(-1)^{k+1}{k+1\choose k-1}\<1>,\\
\<a_{k,k}> =&\ \<k+1>-{k+1\choose 1}\<k>+{k+1\choose 2}\<k-1>-\cdots+(-1)^k{k+1\choose k}\<1> .
\end{align*} 
have the first coordinates equal to Eulerian numbers: $A_{k,i}$ for $i=0,1,\ldots,k-1$ \cite{GKP}.\\
Let us take the cube $C^k=[0,1]^k$ with a normalized  volume equal $k!$.\\
 Then a geometric interpretation of Eulerian number $A_{k,i}$ is a volume 
of i-th slice of the cube $C^k$   located between two successive parallel planes with normal vector [1,1,\ldots,1] \cite{SalKup,EhReSte}. The planes contain consecutive  groups of ${k\choose i}$  vertices of cube  $C^k$. In the geometric interpretation of Eulerian numbers, a cube $C^k$  can be replaced with any k-dimensional parallelepiped with  volume equal $k!$.
So we can interpret expressions  $\<a_{k,k-1}>$ as slices of this k-dimensional parallelepiped.\\

\noindent \noindent Eq. (\ref{14}) is special case of the following equation.
\begin{align}\begin{split}
\<x>= \sum_{j=1}^{k+1}\ \Big(\ \prod_{\substack{i=1\\ i\neq j}}^{k+1}\frac{x-x_i}{x_j-x_i}\  \Big) \<x_j >.
\label{16}
\end{split}\end{align}
\noindent Eq. (\ref{16}) is true in the arithmetic sense but is not true in the geometric sense for any $x,\ x_1,x_2,\ldots,x_{k+1} \in \mathbb{R} $.\\
Probably the following  equation is true in the geometric sense for all $n,m \in \mathbb{Z}$.
\begin{align*}\begin{split}
\<n>= \sum_{j=1}^{k+1}\ \Big(\ \prod_{\substack{i=1\\ i\neq j}}^{k+1}\frac{n-m-i}{j-i}\  \Big) \<m+j >.
\end{split}\end{align*}\\

\section{Multiplication in the set $\mathbb{R}_k$.}

\noindent In this section, we assume that the Eq. (\ref{13}) is true in the geometric sense.
Let us fix a point in $\mathbb{R}_k$ and denote it as $\<0,0,\ldots,0>$. Then, using Eq. (\ref{13}), every element of the set $\mathbb{R}_k$ can be obtained from the distinguished element $\<0,0,\ldots,0>$ as an element $\<x_1,x_2,\ldots,x_{k+1},0>$. We will mark it 
$\<x_1,x_2,\ldots,x_{k+1}>$.\\
For the element $\<t>=\<t_1,t_2,\ldots,t_{k+1},0>=\<t_1,t_2,\ldots,t_{k+1}> $
 Eq. (\ref{13}) takes the following form
 
 \begin{align}\begin{split}
  \forall x_1,x_2,\ldots,x_{k+1},t \in \mathbb{R}  \\
\<x_1,x_2,\ldots,x_{k+1},t>&= \< x_1+t_1,x_2+t_2,\ldots,x_{k+1}+t_{k+1}>\\    
        &=\< x_1+t_1,x_2+t_2,\ldots,x_k+t_k,t_{k+1}> \\
        &\quad  +\dots  + \<t_1, x_2+t_2,x_3+t_3,\ldots,x_{k+1}+t_{k+1}>\\ 
        &\quad - \< x_1+t_1,x_2+t_2,\ldots,x_{k-1}+t_{k-1},t_k,t_{k+1}> \\
     &\quad  -\dots  - \<t_1,t_2, x_3+t_3,x_4+t_4,\ldots,x_{k+1}+t_{k+1}>\\ 
   & \begin{tikzpicture}
         \draw[dotted, line width=1pt] (0,0) -- (8,0);
     \end{tikzpicture} \\
        &\quad  (-1)^{k+1} \big[ \<x_1+t_1,t_2,\ldots,t_{k+1}> \\           
        & \quad + \dots + \< t_1,t_2,\ldots,t_k,x_{k+1}+t_{k+1} > \big] \\       
        &\quad (-1)^k \< t_1,t_2,\ldots,t_{k+1} >.    \label{17}
 \end{split}   \end{align}

\noindent It is easy to see that the elements $\<x_1,x_2,\ldots,x_{k+1}>$ can be interpreted as points of the space $\mathbb{R}^{k+1}$.\\
So all simplices from   $\mathbb{R}_k$ satisfying the condition 
$x_1+x_2+\dots +x_{k+1}=t$ for a fixed $t$, correspond to a hyperplane in
 $\mathbb{R}^{k+1}$.\\
Let us also observe that the elements of Eq. (\ref{17}) form the vertices of a  (k+1)-dimensional rectangular cuboid  in $\mathbb{R}^{k+1}$. \\

\subsection{Multiplication in the set $\mathbb{R}_k$}

Just as in the set $\mathbb{R}_3$, one can also expect that for any natural
number $k > 3$, there exists a multiplication in the set $\mathbb{R}_k$.
By applying the same reasoning from Subsection 3.1 , we obtain the following formulas.\\ 
\[
\<0,\ldots,0,1_i,0\ldots,0>^2= \Big\langle -\frac{1}{k+1},\ldots,
  -\frac{1}{k+1},\frac{2k+1}{k+1},-\frac{1}{k+1},\ldots,-\frac{1}{k+1}\Big\rangle
\] 
\begin{multline*}
\<0,\ldots,0,1_i,0\ldots,0>\<0,\ldots,0,1_j,0\ldots,0>\\
 = \Big\langle -\frac{1}{k+1},\ldots,
  -\frac{1}{k+1},\frac{k}{k+1},-\frac{1}{k+1}\ldots,-\frac{1}{k+1},
  \frac{k}{k+1},-\frac{1}{k+1}\ldots,-\frac{1}{k+1},\Big\rangle,
\end{multline*}
where $1_i, \: 1_j$ denote the number 1 in the i-th  and
 j-th component respectively,\\
  $i\neq j$,\\
number $\frac{2k+1}{k+1}$ are in the in the i-th   component,\\  
 numbers $\frac{k}{k+1}$ are in the in the i-th  and
 j-th component.\\
 
\noindent A general formula for multiplication is as follows
\begin{align}
\<x_1,x_2,\ldots,x_{k+1}>\<a_1,a_2,\ldots,a_{k+1}> =
\<Y_1,Y_2,\ldots,Y_i,\ldots,Y_{k+1}>,
  \label{18}
\end{align}
where
 \begin{multline*}
Y_i
  =   \frac{(-x_1,-\dots -x_{i-1}+(2k+1)x_i-\dots -x_{k+1})(a_1+a_2+\dots     +a_{k+1})}{2(k+1)}\\[0.1cm]
+  \frac{(x_1+x_2+\dots +x_{k+1})(-a_1-\dots -a_{i-1}+(2k+1)a_i-a_{i+1}-\dots -a_{k+1})}{2(k+1)}.
 \end{multline*}
 As in Theorem 3.1, one can show the associativity of the multiplication given by Eq. (\ref{18}).

\end{document}